# Chapter 1:

## Matrix Online Characteristic Number and Goldbach's Conjecture

Baoshan Zhang [1]

(School of Mathematics and Statistics, Nanjing Audit University, Nanjing, China , Zip code: 211815)

Email: bszhang@nau.edu.cn, bszhang8@163.com

**Abstract**：In this paper, it is determined by the consecutive odd numbers between 3 to $2n+1$, and study to the intrinsic properties of a class of matrix sequence $\{A_n\}$. Through the establishment of matrix online number concept, characteristics and the online number column use mathematical induction to prove the some properties of this kind of matrix on the number of online features (**Theorem 1**). Finally, it is given a trial to prove the Goldbach's conjecture (**Theorem 6**). This is the author in the years to explore prime properties in the process of research and discovery, and believe that this finding is of great significance.

**Keywords**　Prime numbers；Online character; Online characteristic columns; Goldbach's conjecture.

**Classification**：O15.156　　**MSC**：11A41

## §1 Introduction

Goldbach 's Conjecture is one of the three Difficult mathematical problems in modern world, It was first found in teaching by C.Goldbach(1690-1764）, who was the German high school teachers [1~3]. On June 7, 1742 , Goldbach write to the great mathematician Euler (1707～1783), and asked that his experiments have found.

Euler said to Goldbach in a reply letter, that he believes that the conclusion is correct, but he can't prove. And then Euler confirmed the following two guess [1~8]:

(1) Strong Goldbach's Conjecture: any more than 4 even can be always expressed as the sum of the two odd primes.

(2) Weak Goldbach's Conjecture: any more than 8 the odd Numbers can always expressed as the sum of the three odd primes.

For many years, this problem caused almost all the attention of mathematicians and others. Goldbach's Conjecture is a mathematical "pearl" in the math crown, which seems a coveted [3].

---

[1] Baoshan Zhang. Male, 1959 born, PhD, Professor, Native of Fengxian County in jiangsu province of China. Mainly engaged in the study of differential equations, computer symbolic operation and so on, also has research interest in number theory, published over 70 papers, published by 2 books. Email: bszhang@nau.edu.cn, bszhang8@163.com

Nearly 300 years passed, and though no Goldbach's conjecture that was completely, but for the theory of the mathematics problem research push, especially the development of analytic theory, get many important results[1,4~18]. In fact, In 1920, Norway's Brown (V.Brun) with screening method to prove: every four even than can be expressed as the sum of the two odd primes, where their factor decomposition number is no more than nine, and it is referred to as "9 + 9". In 1924, H.Rademacher proved "7 + 7". In 1932, T.Estermann proved "6 + 6". In 1937, G.Ricci proved "5 + 7", the "4 + 9", "3 + 15" and "2 + 366". In 1938, A.A.Buchstab proved" 5 + 5 ". In 1940, A.A.Buchstab proved" 4 + 4 ". In 1948, A.Renyi proved "1 + c", which is a great natural number c. In 1956, Wang yuan(China) proved that the "3 + 4". In 1957, Wang yuan(China) proved that the "3 + 3" and "2 + 3". In 1962, Pan Chengdong(China) proved such as "1 + 5," Wang yuan(China) proved "1 + 4". In 1965, A.A.Buchstab proved "1 + 3". In 1966, Chen JingRun(China) proved "1 + 2".

Chen JingRun's work published in 1973 "1 + 2" is the best results, it named as Chen's theorem as below [3, 19].

**Theorem** (Chen's Theorem) Any sufficient big even number is the sum of a prime and a natural number, where the latter is the product of two prime numbers at most.

However, no update progress of Goldbach's conjecture, since "1 + 2" later. Contemporary mathematicians generally thinks that to prove Goldbach's conjecture it must be with the innovation methods, and use new ideas [2, 19~20]. The British "nature" magazine web site reported recently, Prof. Tao Zhexuan (in the university of California , Chinese mathematicians) make out the breakthrough in the study of Goldbach's conjecture weak, and it is expected to finally solve the problem of this century. This is because he has not been in an article in a paper published proof, can be written in an odd number five prime the sum. He said: "the required number of primes is expected to drop to three, it can be proved Goldbach's weak guess. But the Weak Goldbach's guess and compared Goldbach's conjecture or more easily, and so to proof of Goldbach's conjecture, mathematicians maybe to face still huge difficulties."

But just recently, Zhang Yitang gives the weak form of a certificate for the twin prime conjecture, he is undoubtedly the result of the theory, an important milestone in the develop[21~23]. At the same time, weak Goldbach's conjecture has been proved that [24].

This article through the establishment of matrix online characteristic number method，to explore the validity of Goldbach's conjecture. This is the author in the years to explore prime properties in the process of research and discovery, and believe that this finding is of great significance.

## §2 Definitions and Nature of the Theorem

Assume that $N$ express as the natural number set, $P$ as the set of prime, $\Phi$ as an empty set. Let $n \in N$ be an odd natural number. For all $2i+1, 2j+1$ among of $3, 5, 7, 9, \cdots, 2n-1, 2n+1$, Considering the sum $2(i+j+1)$. We costructure matrix $A_n = (a_{ij})_{n \times n}$ with these numbers by

$$a_{ij} = \begin{cases} 2(i+j+1), & \forall 2i+1 \in P \land \forall 2j+1 \in P, \\ 0, & \forall 2i+1 \notin P \lor \forall 2j+1 \notin P. \end{cases} \quad (1)$$

For example, When $n = 8$, the matrix $A_8$ is

$$\begin{pmatrix} 6 & 8 & 10 & 0 & 14 & 16 & 0 & 20 \\ 8 & 10 & 12 & 0 & 16 & 18 & 0 & 22 \\ 10 & 12 & 14 & 0 & 18 & 20 & 0 & 24 \\ 0 & 0 & 0 & 0 & 0 & 0 & 0 & 0 \\ 14 & 16 & 18 & 0 & 22 & 24 & 0 & 28 \\ 16 & 18 & 20 & 0 & 24 & 26 & 0 & 30 \\ 0 & 0 & 0 & 0 & 0 & 0 & 0 & 0 \\ 20 & 22 & 24 & 0 & 28 & 30 & 0 & 34 \end{pmatrix}$$

When $n = 15$, matrix $A_{15}$ is

$$A_{15} = \begin{pmatrix} 6 & 8 & 10 & 0 & 14 & 16 & 0 & 20 & 22 & 0 & 26 & 0 & 0 & 32 & 34 \\ 8 & 10 & 12 & 0 & 16 & 18 & 0 & 22 & 24 & 0 & 28 & 0 & 0 & 34 & 36 \\ 10 & 12 & 14 & 0 & 18 & 20 & 0 & 24 & 26 & 0 & 30 & 0 & 0 & 36 & 38 \\ 0 & 0 & 0 & 0 & 0 & 0 & 0 & 0 & 0 & 0 & 0 & 0 & 0 & 0 & 0 \\ 14 & 16 & 18 & 0 & 22 & 24 & 0 & 28 & 30 & 0 & 34 & 0 & 0 & 40 & 42 \\ 16 & 18 & 20 & 0 & 24 & 26 & 0 & 30 & 32 & 0 & 36 & 0 & 0 & 42 & 44 \\ 0 & 0 & 0 & 0 & 0 & 0 & 0 & 0 & 0 & 0 & 0 & 0 & 0 & 0 & 0 \\ 20 & 22 & 24 & 0 & 28 & 30 & 0 & 34 & 36 & 0 & 40 & 0 & 0 & 46 & 48 \\ 22 & 24 & 26 & 0 & 30 & 32 & 0 & 36 & 38 & 0 & 42 & 0 & 0 & 48 & 50 \\ 0 & 0 & 0 & 0 & 0 & 0 & 0 & 0 & 0 & 0 & 0 & 0 & 0 & 0 & 0 \\ 26 & 28 & 30 & 0 & 34 & 36 & 0 & 40 & 42 & 0 & 46 & 0 & 0 & 52 & 54 \\ 0 & 0 & 0 & 0 & 0 & 0 & 0 & 0 & 0 & 0 & 0 & 0 & 0 & 0 & 0 \\ 0 & 0 & 0 & 0 & 0 & 0 & 0 & 0 & 0 & 0 & 0 & 0 & 0 & 0 & 0 \\ 32 & 34 & 36 & 0 & 40 & 42 & 0 & 46 & 48 & 0 & 52 & 0 & 0 & 58 & 60 \\ 34 & 36 & 38 & 0 & 42 & 44 & 0 & 48 & 50 & 0 & 54 & 0 & 0 & 60 & 62 \end{pmatrix} \quad (2)$$

The matrix (2) has the following distinctive features:

（1）$A_{15}$ is a symmetric matrix, and the position $i, j$ of 0 in $A_{15}$ being satisfy that $2i+1, 2j+1$ at least one is not a prime number.

（2）Each is an even number in $A_{15}$ is equal to the location $i, j$ that corresponds to the sum of two prime numbers $2i+1, 2j+1$.

（3）Each even number in $A_{15}$ has the properties of the same even in with pair of parallel to the deputy diagonal line.

（4）Except 2,3,56, all even numbers less than or equal to 62 appeared in $A_{15}$, they grew up in with pair of parallel to the deputy diagonal line in turn.

For being easy to describe, we define the concept of the matrix online digital features.

**Definition 1** For the $n \times n$ matrix $A$, the count of all null number it with the same pair of diagonal parallel each online is called **online digital characteristics**. Starting from the top left corner to the lower right corner of the end of online digital characteristics of matrix of array called **online digital characteristic columns** of $A$.

Obviously, the online digital characteristic is always a negative integer, and it is equal to zero if and only if all digital 0 in this online.

Let $L_A(k)$ be the $kth$ online digital characteristic of the matrix $A$, which is from the upper left corner to the lower right corner. Then, the online digital characteristics of the matrix $A$ are consists of a number of $2n-1$, respectively by

$$L_A(1), \ L_A(2), \cdots, L_A(n), \ L_A(n+1), \cdots, L_A(2n-1). \qquad (3)$$

For example, for given above $A_{15}$, the online digital characteristic columns of $A_{15}$ is that

$$1, 2, 3, 2, 3, 4, 4, 4, 5, 6, 5, 4, 6, 4, 7, 8, 3, 4, 6, 2, 3, 4, 2, 2, 2, 0, 1, 2, 1 \qquad (4)$$

It worth pointing out that, the online digital characteristic of $A_{15}$ reflect that even number referred to in the online can be expressed as the sum of two prime numbers, this is the number of possible pairs of prime. Specifically, $L_{A_{15}}(k)$, that is the $kth$ number in (4), represents an even number $2(k+2)$ can be expressed as the sum of two prime numbers of $L_{A_{15}}(k)$ primes.

For instance, $L_{A_{15}}(16) = 8$ means $2(16+2) = 36$ can be expressed as the sum of two prime numbers of eghit primes. In fact, we have

$$36 = 5 + 31, \ 36 = 7 + 29, \ 36 = 13 + 23, \ 36 = 17 + 19,$$
$$36 = 31 + 5, \ 36 = 29 + 7, \ 36 = 23 + 13, \ 36 = 19 + 17.$$

It is easy to know, that the number of online features columns (4) in $A_{15}$, there is only one zero, and its position of appeared in the later. To describe this phenomenon, we put the line number feature column is divided into fluctuation two parts (according to the deputy diagonal lines are divided into two groups), and now we introduce the following definition.

**Definition 2** For the online digital characteristic columns of $n \times n$ Matrix $A$, the former $n$ numbers

$$L_A(1), \ L_A(2), \cdots, L_A(n-1), L_A(n) \quad (5)$$

is called **the upper left corner online digital characteristic columns** of $A$. And the later $n-1$ numbers

$$L_A(n+1), L_A(n+2), \cdots, L_A(2n-1) \quad (6)$$

is called **the lower right corner online digital characteristic columns** of $A$.

From formula (4), we obtain, the upper left corner online digital characteristic columns of $A_{15}$ is

1, 2, 3, 2, 3, 4, 4, 4, 5, 6, 5, 4, 6, 4, 7.

And the lower right corner online digital characteristic columns is

8, 3, 4, 6, 2, 3, 4, 2, 2, 2, 0, 1, 2, 1.

It is not difficult to find that, the upper left corner of the line number feature of $A_{15}$ are all greater than zero, and its lower right corner all online digital characteristics are almost greater than zero. Meanwhile, the zero appear behind in this columns (near the position of $L_{A_{15}}(29)$). To Described here, we asked for general defined by (1) the characteristics of the matrix $A_n$ as well? Let us discuss the problem as below.

Easy to verify that, for $n = 1, 2, 3, \cdots, 14$, the corresponding matrix $A_1, A_2, A_3, \cdots, A_{14}$ have the nature of the number of online characteristics form as $A_{15}$ we discussed. Therefore, we want to prove that, for all the shape, such as definition (1) of the matrix $A_n$ have such properties. That is to prove the following theorem.

**Theorem 1** For any natural number $n \in N$, the matrix $A_n$, such as definition (1), has always the nature of below: The upper left corner number online digital characteristics of $A_n$ are all greater than zero. and its lower right corner number online digital characteristics of $A_n$ is greater than zero almost, Meanwhile, the zero appear behind in this columns $\{L_{A_n}(j)\}_{j=1}^{2n-1}$ (near the position of $L_{A_n}(2n-1)$).

To prove this theorem, we first discuss the intrinsic characteristics of the matrix sequence $\{A_n\}$ defined by (1).

Let us calculate the first 4 items of the matrix sequence $\{A_n\}$:

$$A_1 = (6), \quad A_2 = \begin{pmatrix} 6 & 8 \\ 8 & 10 \end{pmatrix}, \quad A_3 = \begin{pmatrix} 6 & 8 & 10 \\ 8 & 10 & 12 \\ 10 & 12 & 14 \end{pmatrix}, \quad A_4 = \begin{pmatrix} 6 & 8 & 10 & 0 \\ 8 & 10 & 12 & 0 \\ 10 & 12 & 14 & 0 \\ 0 & 0 & 0 & 0 \end{pmatrix}. \quad (7)$$

We can see that, $A_1, A_2, A_3$ are the ordinal master submatrix of $A_4$, and $A_1, A_2, A_3, A_4$ are also the ordinal master submatrix of $A_{15}$. It is evident that this result has a general. So we have

**Theorem 2** For any natural number $n \in N$, the matrix $A_n$, such as definition (1), has always the peculiarity: that $A_1, A_2, A_3, \cdots, A_{n-1}$ are the ordinal master submatrix of $A_n$.

Using this result, we proved that

**Theorem 3** For any natural number $n \in N$, the matrix $A_n$, such as definition (1), has always the peculiarity: that the former $n-1$ numbers of the top left corner of the online digital characteristic columns of $A_n$ is consist of the top left corner of the online digital characteristic columns of $A_{n-1}$.

**Proof** According to (5), the top left corner of the online digital characteristic columns of $A_n$ is

$$L_{A_n}(1), \quad L_{A_n}(2), \cdots, L_{A_n}(n-1), L_{A_n}(n).$$

And the top left corner of the online digital characteristic columns of $A_{n-1}$ is

$$L_{A_{n-1}}(1), \quad L_{A_{n-1}}(2), \cdots, L_{A_{n-1}}(n-1).$$

Using theorem 2, we have

$$L_{A_n}(j) = L_{A_{n-1}}(j), \quad j = 1, 2, \cdots, n-1. \quad (8)$$

This is the conclusion of theorem 3.

**Theorem 4** For any natural number $n \in N$, the matrix $A_n$, such as definition (1), has always the peculiarity: The former $n-3$ numbers of the lower right corner of the online digital characteristic columns of $A_n$ is successively not less than the later $n-3$ numbers of the lower right corner of the online digital characteristic columns of $A_{n-1}$. At the same time, the online diagonal characteristic $L_{A_n}(n)$ of $A_n$ is successively not less than $L_{A_{n-1}}(n)$.

**Proof** Since the online diagonal characteristics consist of a number of $2n-1$ for a $n \times n$ matrix, and the online diagonal characteristics consist of a number of $2n-3$ for a $(n-1) \times (n-1)$ matrix. According to (5), the lower right corner of the online digital characteristic columns of $A_n$ is

$$L_A(n+1), L_A(n+2), \cdots, L_A(2n-1).$$

And the lower right corner of the online digital characteristic columns of $A_{n-1}$ is

$$L_{A_{n-1}}(n), L_{A_{n-1}}(n+1), \cdots, L_{A_{n-1}}(2n-3).$$

Using theorem 2, and the definition of the online digital characteristic, we have

$$L_{A_n}(j) \geq L_{A_{n-1}}(j), \quad j = n, n+1, n+2, \cdots, 2n-3. \tag{9}$$

This is the conclusion of theorem 4.

Composite of theorem 3, theorem 4, that is formula (8), (9), we get the following conclusion.

**Theorem 5** For any natural number $n \in N$, the matrix $A_{n-1}, A_n$, such as definition (1), the online diagonal characteristics of $A_n$ satisfy

$$L_{A_n}(j) = \begin{cases} L_{A_{n-1}}(j), & j = 1, 2, \cdots, n-1, \\ L_{A_{n-1}}(j) \text{ or } L_{A_{n-1}}(j) + 2, & j = n, n+1, n+2, \cdots, 2n-3, \\ 0 \text{ or } 2, & j = 2n-2, \\ 0 \text{ or } 1, & j = 2n-1. \end{cases} \tag{10}$$

**Proof** Since the order of Matrix $A_n$ is bigger 1 than the order of Matrix $A_{n-1}$, the corresponding online digital characteristic columns increase 2. The lower right corner of online characteristic of $A_n$ changes depends on the digital values about the columns and the lines of $A_n$. According to the symmetry of $A_n$, for $n \leq j \leq 2n-3$, based on $L_{A_{n-1}}(j)$, $L_{A_n}(j)$ either not increase (especially when the $2n+1$ is composite!) or increase 2 (at this time, $2n+1$ must be prime !). Thus, according to (8), (9), we get the two expressions of (10).

For $2n-2 \leq j \leq 2n-1$, the value of $L_{A_n}(j)$ depends on whether the $2n-1, 2n+1$ is the prime situation.

If $2n-1, 2n+1$ are all prime, then $L_{A_n}(2n-2) = 2$, $L_{A_n}(2n-1) = 1$.

If $2n-1$ is composite and $2n+1$ is prime, then $L_{A_n}(2n-2) = 0$, $L_{A_n}(2n-1) = 1$.

If $2n+1$ is composite, then $L_{A_n}(2n-2) = 0$, $L_{A_n}(2n-1) = 0$.

Comprehensive the above discussion, we have acquired the proof of the theorem 5.

Theorem 5 depicts the relation of two adjacent matrix of the sequence $\{A_n\}$ defined by (1) about their online digital characteristic columns. We will see, these results will help us to complete the proof of theorem 1.

### §3 Proof of theorem 1 and its application

Now consider the proof of theorem 1.

**Proof of Theorem 1** For the matrix sequence $\{A_n\}$ of Defined by (1), We use mathematical induction to $n \in N$.

First of all, when $n \in N$ is small, the proposition is true. In fact, by (7) can be verified directly. That is for small $n \in N$, the upper left corner number online digital characteristics of $A_n$ are all greater than zero. and its lower right corner number online digital characteristics of $A_n$ is greater than zero almost, Meanwhile, the zero appear behind in this columns $\{L_{A_n}(j)\}_{j=1}^{2n-1}$ (near the position of $L_{A_n}(2n-1)$). Especially, when $n = 15$, we know that proposition is true by (4).

Secondly, the induction hypothesis, If $n = m$, the proposition is true. That is for $m \in N$, the upper left corner number online digital characteristics of $A_m$ are all greater than zero. and the lower right corner number online digital characteristics of $A_m$ are greater than zero almost, Meanwhile, the zero appear behind in this columns $\{L_{A_m}(j)\}_{j=1}^{2m-1}$ (near the position of $L_{A_m}(2m-1)$).

Here, the lower right corner number online digital characteristics of $A_m$ is greater than zero almost, that means the number $m_0$ of zero in $\{L_{A_m}(j)\}_{j=m+1}^{2m-1}$ is trivial compared with $m$, thus $\frac{m_0}{m} \approx 0$. But for the lower right corner number online digital characteristics columns of $A_m$, the zreo position $k_0$ which its appear first time in $\{L_{A_m}(j)\}_{j=m+1}^{2m-1}$, is near the number $2m-1$, that is $\dfrac{k_0}{2m-1} \approx 1$.

Now considering the case $n = m+1$, Whether the proposition is true or not. For the matrix $A_{m+1}$ defined by (1), since $A_m$ is the *mth* ordinal master submatrix of $A_{m+1}$, and $A_m$ is a $m \times m$ matrix defined by (1), and so it conforms to the induction hypothesis. Thus we have

$$L_{A_m}(j) > 0, \quad j = 1, 2, \cdots, m-1, m, m+1. \tag{11}$$

Therefore by the formula (10) in theorem 5, we obtain that the upper left corner number online digital characteristics of $A_{m+1}$ will satisfy

$$L_{A_{m+1}}(j) = L_{A_m}(j) > 0 \; (j = 1, 2, \cdots, m-1, m), \quad L_{A_{m+1}}(m+1) \geq L_{A_m}(m+1) > 0. \tag{12}$$

And so the upper left corner number online digital characteristics of $A_{m+1}$ are all more than zero.

In another hand, since the lower right corner number online digital characteristics columns of $A_m$ are almost more than zero, thus the lower right corner number online digital characteristics columns of $A_{m+1}$ are almost more than zero too. Moreover, Let $\tilde{m}_0$ be the number of zero in $\{L_{A_{m+1}}(j)\}_{j=m+2}^{2m+1}$, $\tilde{k}_0$ be the zreo position which its appear first time in $\{L_{A_{m+1}}(j)\}_{j=m+2}^{2m+1}$.

According to the formula (10) in theorem 5 and its prove that reveals inner link between the lower right corner number online digital characteristics columns $A_{m+1}$ and $A_m$, we have

$$0 \leq \tilde{m}_0 \leq m_0 + 2, \quad 2m+1 \geq \tilde{k}_0 \geq k_0. \tag{13}$$

Here, $m_0$ be the number of zero in $\{L_{A_m}(j)\}_{j=m+1}^{2m-1}$, $k_0$ be the zreo position which its appear first time in $\{L_{A_m}(j)\}_{j=m+1}^{2m-1}$. So according to the induction hypothesis, it would had $\dfrac{m_0}{m} \approx 0$ and $\dfrac{k_0}{2m-1} \approx 1$. Thus we have

$$0 \leq \dfrac{\tilde{m}_0}{m+1} \leq \dfrac{m_0 + 2}{m+1} = \dfrac{m_0}{m+1} + \dfrac{2}{m+1} \approx \dfrac{2}{m+1}, \tag{14}$$

$$1 \geq \dfrac{\tilde{k}_0}{2m+1} \geq \dfrac{k_0}{2m+1} \approx \dfrac{k_0}{2m-1} \approx 1. \tag{15}$$

Here might as well set $m$ is large enough (for instance $m$ more than $10^5$), By (14), (15), we have

$$\frac{\tilde{m}_0}{m+1} \approx 0, \quad \frac{\tilde{k}_0}{2m+1} \approx 1. \quad (16)$$

This shows that $\tilde{m}_0$ is insignificant in comparison with $m+1$, and $\tilde{k}_0$, $2m+1$ are very close. Therefore the proposition t is true for $n = m+1$.

Finally, according to the principle of mathematical induction, theorem 1 is true for any natural number $n \in N$.

The conclusion of theorem 1 is very important. Because for any natural number $n \in N$, the upper left corner number online digital characteristics of $A_n$ defined by (1) are all greater than zero, and so, we can know that any even number $2m$ between 6 and $2n+4$ can be expressed as the sum of two prime numbers $p_1, p_2$, where $p_1, p_2$ come from between 3 and $2m$. That is to say, $2m = p_1 + p_2$. Specially, $2n = p + q$, where $p, q$ are prime between 3 and $2n$. This is the result of Goldbach's conjecture to be prove that.

As the application of theorem 1, we now put the Goldbach's conjecture, account for the following theorem and its proof is given.

**Theorem 6** (**Goldbach's conjecture**) Any even number larger than 4 can be expressed as the sum of two odd prime numbers. That is, for any natural number $n(\geq 3) \in N$, there is always a pair of prime numbers $p_1, p_2 \in P$, such as $2n = p_1 + p_2$.

**Proof** For any natural number $n(\geq 3) \in N$, it might assume $n$ is large enough. According to the formula (1), we construct $(n-2) \times (n-2)$ matrix $A_{n-2}$.

By theorem 1, the upper left corner number online digital characteristics of $A_{n-2}$ are all greater than zero. Of course, it must be $L_{A_{n-2}}(n-2) > 0$. That is to say, number $2(n-2)+4 = 2n$ appear at least once on the vice diagonal line of $A_{n-2}$.

In this case, there is pair of prime numbers $p_1, p_2 \in P$ between 3 and $2(n-2)+1 = 2n-3$, that is satisfying $2n = p_1 + p_2$.

This is a simple proof of theorem 6.

## §4 Conclusions and Remarks

The proof of theorem 6 is based on theorem 1. However, the proof of theorem 1 is done by using mathematical induction. In these discussions, in the case of depicting two primes and said even we set up concept of matrix online digital characteristics and online digital characteristics column, it brought convenience to discuss the problem. If theorem 1 proved to be correct, the proof of theorem 6 is the Goldbach's conjecture to prove, that is to say, the puzzled people for a long time of number theory problem was solved. This, of course, is a very happy thing.

By the way, the twin prime conjecture is also has yet to solve arithmetic problems in the world now.

**Conjecture 1**[23,24]（Twin prime conjecture）There are infinitely many prime numbers $p$, $p$ with $p+2$ is a prime number; And prime number $(p, p+2)$ called the twin prime number.

Mathematicians believed that the above conjecture is true, but it has not been proved completely. Zhang Yitang's job today is to thoroughly solve the weak twin prime conjecture made outstanding contributions. The author tries to consider the sequence sets

$$S_a = \{ k \in N \mid 2k + a \in P \},  \quad (17)$$

where $a \in N$ is any odd number, and by using similar methods to prove that $S_a \cap S_{a+2}$ is not for an empty set in my other thesis. And eventually, we try to make out the proof of the twin prime conjecture. This work, of course, might to be the mathematics examination and judgment of authority.

Further, if always prove for $S_a \cap S_{a+2}$ being the empty set is correct, it seems to promote more generally, Polignac's conjecture [23,24] (Alphonsede Polignac, 1817 ~ 1890)：

**Conjecture 2**（Polignac's conjecture）For all natural number $k$, there are infinitely many prime numbers of $(p, p+2k)$.

In fact, the author seems to confirm that for any given number $k$, as long as there is a pair of primes $(p_0, p_0 + 2k)$, it must prove there are infinitely many primes $(p, p+2k)$. Of course, this work remains to be further discussed, so dare not presume to here.

# Chapter 2:

# Matrix Master Characteristic Sequence and Polignac Conjecture

**Baoshan Zhang**[①]

(School of Mathematics and Statistics, Nanjing Audit University, Nanjing, China, Zip code: 211815)

Email: bszhang@nau.edu.cn  bszhang8@163.com

**Abstract**: In this paper, it is defined the concepts of matrix master characteristic number and the Matrix Master Characteristic Sequence (**Definition 1**). Firstly, we prove that any even number can be expressed as for the difference of two odd prime numbers at least two groups (**Theorem 4**). Secondly, we prove that there are infinitely many odd prime numbers separated by four (**Theorem 9**). Finally, we prove that if there is greater than 1 in $S_3 \cap S_{2m+3}$ for any $m \in N$, so that there are infinitely many odd prime numbers separated by $2m$ (**Theorem 11**). The results will undoubtedly promote the question's research for Polignac conjecture.

**Keywords:** Prime numbers; The twin prime conjecture; Matrix master characteristic number; Matrix master characteristic sequence; Polignac conjecture

**MSC**: 11A41

## §1  Introduction

It is well known that there are two famous conjecture about prime problems as below [1~6]。

**Conjecture 1** (**Twin prime conjecture**) there are infinitely many prime numbers $p$, such as $p$ and $p+2$ are primes. And the pair prime $(p, p+2)$ called the twin prime number.

**Conjecture 2** (**Polignac conjecture**) For all natural numbers $k$, there are infinitely many prime numbers $p$, such as $p$ and $p+2k$ are prime.

Mathematicians believed that the above conjectures are all true, but it has not been proved completely. On May 14, 2013, Natural Magazine reported online Chinese mathematician Zhang Yitang proves that "there are an infinite number of primes less than 70 million ". This research then is considered on the twin prime conjecture the ultimate number theory problem made a major breakthrough, even some people think that their impact on the community will be more than the "1 + 2" trained Chen jingrun proved [6~23]. Although the Zhang Yitang is to thoroughly solve the twin prime conjecture made outstanding contributions, but ultimately solve the twin prime conjecture is still a considerable distance.

That is to say, the twin prime conjecture is still not yet thoroughly crack number theory problem, the Polignac conjecture is also can't be solved [21~24]。 Recently, the author in another

---

paper, consider the set

$$S_a = \{ k \in N \mid 2k + a \in P \}, \tag{1}$$

where $N$ is natural set, $P$ is prime set, and $a \in N$ is any odd number. And try to prove the set $S_a \cap S_{a+2}$ is not empty, and then to solve the twin prime conjecture. However, this work needs to be international math examination and judgment of authority. This article is based on the above work continues but from another perspective, and make some discusses. The purpose, of course, is to promote finding the answer of Polignac conjecture.

## §2 Definition and Preparatory Theorems

For a given natural number $n$, let consider the absolute value of difference $2|i-j|$ of the odd numbers $2i+1, 2j+1$ between 3 to $2n+1$, and structure matrix $A_n = (a_{ij})_{n \times n}$, where

$$a_{ij} = \begin{cases} 2|i-j|, & \forall 2i+1 \in P \land \forall 2j+1 \in P, \\ 0, & \forall 2i+1 \notin P \lor \forall 2j+1 \notin P. \end{cases} \tag{2}$$

For example, when $n = 15$, the matrix $A_{15}$ is

$$A_{15} = \begin{pmatrix}
0 & 2 & 4 & 0 & 8 & 10 & 0 & 14 & 16 & 0 & 20 & 0 & 0 & 26 & 28 \\
2 & 0 & 2 & 0 & 6 & 8 & 0 & 12 & 14 & 0 & 18 & 0 & 0 & 24 & 26 \\
4 & 2 & 0 & 0 & 4 & 6 & 0 & 10 & 12 & 0 & 16 & 0 & 0 & 22 & 24 \\
0 & 0 & 0 & 0 & 0 & 0 & 0 & 0 & 0 & 0 & 0 & 0 & 0 & 0 & 0 \\
8 & 6 & 4 & 0 & 0 & 2 & 0 & 6 & 8 & 0 & 12 & 0 & 0 & 18 & 20 \\
10 & 8 & 6 & 0 & 2 & 0 & 0 & 4 & 6 & 0 & 10 & 0 & 0 & 16 & 18 \\
0 & 0 & 0 & 0 & 0 & 0 & 0 & 0 & 0 & 0 & 0 & 0 & 0 & 0 & 0 \\
14 & 12 & 10 & 0 & 6 & 4 & 0 & 0 & 2 & 0 & 6 & 0 & 0 & 12 & 14 \\
16 & 14 & 12 & 0 & 8 & 6 & 0 & 2 & 0 & 0 & 4 & 0 & 0 & 10 & 12 \\
0 & 0 & 0 & 0 & 0 & 0 & 0 & 0 & 0 & 0 & 0 & 0 & 0 & 0 & 0 \\
20 & 18 & 16 & 0 & 12 & 10 & 0 & 6 & 4 & 0 & 0 & 0 & 0 & 6 & 8 \\
0 & 0 & 0 & 0 & 0 & 0 & 0 & 0 & 0 & 0 & 0 & 0 & 0 & 0 & 0 \\
0 & 0 & 0 & 0 & 0 & 0 & 0 & 0 & 0 & 0 & 0 & 0 & 0 & 0 & 0 \\
26 & 24 & 22 & 0 & 18 & 16 & 0 & 12 & 10 & 0 & 6 & 0 & 0 & 0 & 2 \\
28 & 26 & 26 & 0 & 20 & 18 & 0 & 14 & 12 & 0 & 8 & 0 & 0 & 2 & 0
\end{pmatrix} \tag{3}$$

This matrix has the following distinctive features:

(1) $A_{15}$ is a symmetric matrix, any one is zero on the diagonal, and the zero in location $i, j$ meet with $2i+1, 2j+1$ have at least one is not a prime number.

(2) At the position $i, j$, each positive even number of $A_{15}$ is equal to the absolute value of

the difference between the two prime numbers $2i+1, 2j+1$. That is $2|i-j|$.

(3) All same even number are location of the gradient line which parallel to the diagonal line.

(4) All even numbers less than or equal to 28 appear in $A_{15}$, they, based on diagonal, grew up in turn to the diagonal line on both sides of the diffusion.

For easy to describe, we introduced the concept of matrix master characteristic number and master characteristic sequence.

**Definition 1** For $n \times n$ matrix $A$, the count of positive even number, at the same parallel to the diagonal of each line, called its **matrix master characteristic number** of that line. These characteristic numbers array from diagonal to the top right corner called **matrix master characteristic number sequence**.

According to this definition, all matrix characteristic numbers are negative integer. If and only if a certain characteristics number of $A$ is zero, all to zero corresponding to the numbers on that line.

Let $f_A(k)$ be the $k$ th matrix master characteristic number of $A$ from diagonal to the top right corner, then the master characteristic sequence of $A$ consists of $n$ numbers, respectively by

$$f_A(1), \ f_A(2), \cdots, f_A(k), \cdots, f_A(n-1), \ f_A(n). \tag{4}$$

To above matrix $A_{15}$, for instance, the master characteristic sequence of $A_{15}$ is

$$0, 5, 4, 6, 4, 4, 5, 3, 3, 3, 2, 1, 2, 2, 1. \tag{5}$$

It worth pointing out that, the matrix master characteristic number of $A_{15}$ reflects how many the possible pair prime numbers, which the even number can be expressed as the difference referred to in the online two primes (big prime number decrease small prime number, similarly hereinafter).

In particular, $f_{A_{15}}(k)$ means the even number $2(k-1)$ can be expressed as the difference two primes total of $f_{A_{15}}(k)$. By formula (5), there is $f_{A_{15}}(7) = 5$, thus $2(7-1) = 12$ can be expressed as the difference two primes total of 5. In fact, we have

$$12 = 17 - 5, \ 12 = 19 - 7, \ 12 = 23 - 11, \ 12 = 29 - 17, \ 12 = 31 - 19.$$

It is note that $f_{A_{15}}(1) = 0$, This is due in accordance with the determined (2) by the construction of the type defined matrix $A_n = (a_{ij})_{n \times n}$, it is general. Therefore, we only discuss the sequence $\{f_{A_n}(k)\}_{k=2}^{n}$ to the matrix master characteristics of $A$ defined (4). For matrix $A_{15}$, $\{f_{A_{15}}(k)\}_{k=2}^{15}$ are all positive, but it's not general. Because consider $A_{13}$, that is the order master submatrix of matrix $A_{15}$, it has $\{f_{A_{13}}(k)\}_{k=2}^{13}$ as below.

$$4, 4, 5, 3, 3, 3, 2, 2, 1, 1, 0, 0. \qquad (6)$$

Here $f_{A_{13}}(12) = f_{A_{13}}(13) = 0$, and so $\{f_{A_{13}}(k)\}_{k=2}^{13}$ are not always all positive.

In spite of this, it is not hard to find, the number "1, 0" appeared in $\{f_{A_{13}}(k)\}_{k=2}^{13}$ are location of the later position, and the number of "1, 0" is relatively small. In fact, there is only 4 about "1, 0" in $\{f_{A_{13}}(k)\}_{k=2}^{13}$, others all positive, and the first location of "1 or 0" is $k_0 = 9$. Let consider the frequency $\mu$ for "1, 0" and the degree $\nu$ of close to the end of the sequence:

$$\mu = \frac{4}{12} = \frac{1}{3}, \quad \nu = \frac{k_0}{12} = \frac{9}{12} = \frac{3}{4}.$$

Therefore, it hast $\mu \approx 0, \nu \approx 1$. And so we describe the phenomenon $\{f_{A_{13}}(k)\}_{k=2}^{13}$ is almost all greater than 1, and "1, 0" are almost close to the end of $\{f_{A_{13}}(k)\}_{k=2}^{13}$.

Generally, we introduce the following definition.

**Definition 2** For $n \times n$ matrix $A$, Assume $\alpha_n$ be the number of "1, 0" appeared in the master characteristic subsequence $\{f_{A_n}(k)\}_{k=2}^{n}$ of $A$, and $t_n$ be the first location of "1 or 0" appeared in $\{f_{A_n}(k)\}_{k=2}^{n}$. If there are

$$\mu = \frac{\alpha_n}{n-1} \approx 0, \quad \nu = \frac{t_n}{n-1} \approx 1. \qquad (7)$$

then it called the matrix master characteristic number of $A$ is **the almost greater than** 1, and "1, 0" are **almost close to the end** of $\{f_{A_n}(k)\}_{k=2}^{n}$.

Now consider the general matrix sequence $\{A_n\}_{n=2}^{\infty}$ defined by (2), and discuss whether each matrix master characteristic number is the almost greater than 1, and "1, 0" are almost close to the end or not.

Based on the concept of order master submatrix, it is easy to find that matrix sequence $\{A_n\}_{n=2}^{\infty}$ has the following properties:

**Theorem 1** For nay $n \in N$, The former matrix is the order master matrix of the later matrix in the matrices $A_1, A_2, A_3, \cdots, A_{n-1}, A_n$ defined by (2), and all the order master matrix of $A_n$ Followed by $A_1, A_2, A_3, \cdots, A_{n-1}, A_n$.

**Proof** It gets immediately by the matrix structure and the concept of order master Matrix sequence.

For the matrix sequence $\{A_n\}_{n=2}^{\infty}$ defined by (2), we will set up and prove finally the important conclusion as below.

**Theorem 2** For any $n(\geq 2) \in N$, the matrix $A_n$ defined by (2) has the quality: Its matrix master characteristic subsequence $\{f_{A_n}(k)\}_{k=2}^{n}$ is almost greater than 1, and "1, 0" are almost close to the end.

To prove this theorem, we first prove the following auxiliary theorem.

**Theorem 3** For any $n(\geq 2) \in N$, the matrix $A_n$ defined by (2) has the quality: Its matrix master characteristic subsequence $\{f_{A_n}(k)\}_{k=2}^{n}$ satisfies the following relation with the matrix master characteristic subsequence $\{f_{A_{n-1}}(k)\}_{k=2}^{n-1}$ of $A_n$'s order master matrix $A_{n-1}$.

$$f_{A_n}(k) \geq f_{A_{n-1}}(k) \quad (k=2,3,\cdots n-1), \quad f_{A_n}(n) \geq 0. \tag{8}$$

**Proof** Since $A_{n-1}$ is $(n-1)th$ order master matrix of $A_n$, by definition 1, $\{f_{A_n}(k)\}_{k=2}^{n}$ has the following relation with $\{f_{A_{n-1}}(k)\}_{k=2}^{n-1}$:

$$f_{A_n}(k) = f_{A_{n-1}}(k) + sign(a_{n-k+1,n}), k=2,3,\cdots n-1, \quad f_{A_n}(n) = sign(a_{1n}). \tag{9}$$

Here $(a_{1n}, a_{2n}, \cdots, a_{nn})^T$ is the $n$ th column of $A_n$, $sign$ is symbolic function. And so $(a_{1n}, a_{2n}, \cdots, a_{nn})^T \geq 0$, that is all elements of this column are nonnegative.

By formula (2), if $2n+1$ is a prime, then $a_{1n}, a_{2n}, \cdots, a_{n-1n}$ are not all zero, expressly $a_{1n} = 2(n-1)$. If $2n+1$ is a composite, then $a_{1n}, a_{2n}, \cdots, a_{n-1n}$ are all zero at nonce. In any case, the values of symbolic function $sign(a_{1n}), sign(a_{2n}), \cdots, sign(a_{n-1n})$ are all nonnegative. Therefore, we obtain formula (8) immediately by formula (9). This is the conclusion of theorem 3.

In fact, From the process of above proof, we can put the refinement to (8) as below.

$$f_{A_n}(k) = \begin{cases} f_{A_{n-1}}(k), & 2n+1 \notin P \\ \geq f_{A_{n-1}}(k), & 2n+1 \in P. \end{cases} (2 \leq k \leq n-1), \quad f_{A_n}(n) = \begin{cases} 0, & 2n+1 \notin P \\ 1, & 2n+1 \in P. \end{cases} \tag{10}$$

where $P$ is the prime set.

Now we give the proof of theorem 2.

**Proof of Theorem 2** For the matrix sequence $\{A_n\}$ defined by (2), we use mathematical induction for $n \in N$.

Firstly, since above discusses to $A_{13}$ (see formula (6)), we have know for small $n \in N$, the proposition of theorem 2 is true.

Secondly, Let's make the induction hypothesis: Assume that the result of theorem 2 is true when $n = m$. That is, for the matrix $A_m$ defined by (2), its matrix master characteristic subsequence $\{f_{A_m}(j)\}_{j=2}^{m}$ is almost grater than 1, and "1, 0" are almost close to the end $f_{A_m}(m)$.

By definition 2, let $m_0$ be the number of "1, 0" appear in $\{f_{A_m}(j)\}_{j=2}^{m}$, and $k_0$ be the first position of "1 or 0" in subsequence $\{f_{A_m}(j)\}_{j=2}^{m}$. Thus $m_0$ is almost grater than 1, and $k_0$ almost close to $m$. By describe in mathematical language, there will be

$$\mu_m = \frac{m_0}{m} \approx 0, \quad v_m = \frac{k_0}{m} \approx 1.$$

Now let's discuss whether the proposition of theorem 2 is true for $n = m+1$ or not.

For the $(m+1) \times (m+1)$ matrix $A_{m+1}$ defined by (2), we consider its $m$ th order master matrix $A_m$. Since $A_m$ is also a $m \times m$ matrix defined by (2), it will satisfy the induction hypothesis. And so according to the theorem 3, there is

$$f_{A_{m+1}}(j) \geq f_{A_m}(j), (j = 2, \cdots, m-1, m), \quad f_{A_{m+1}}(m+1) \geq 0. \quad (11)$$

According to the induction hypothesis: $\{f_{A_m}(j)\}_{j=2}^{m}$ is almost greater than 1, and "1, 0" are almost close to the end. Therefore, from above mention and formula (11), $\{f_{A_{m+1}}(j)\}_{j=2}^{m+1}$ is almost greater than 1, and "1, 0" are almost close to the end too.

In fact, let $\tilde{m}_0$ be the number of "1, 0" appear in $\{f_{A_{m+1}}(j)\}_{j=2}^{m+1}$, and $\tilde{k}_0$ be the first position of "1 or 0" in $\{f_{A_{m+1}}(j)\}_{j=2}^{m+1}$. According to the formula (10) in theorem 3, there is

$$0 \leq \tilde{m}_0 \leq m_0 + 1, \quad m+1 \geq \tilde{k}_0 \geq k_0. \quad (12)$$

Where $m_0$ is the number of "1, 0" appear in $\{f_{A_m}(j)\}_{j=2}^{m}$, and $k_0$ is the first position of "1 or 0" in subsequence $\{f_{A_m}(j)\}_{j=2}^{m}$. Then according to the induction hypothesis and defined 2, there is

$$0 \leq \mu_{m+1} = \frac{\tilde{m}_0}{m+1} \leq \frac{m_0+1}{m+1} = \frac{m_0}{m+1} + \frac{1}{m+1} = \frac{m}{m+1}\mu_m + \frac{1}{m+1}, \quad (13)$$

$$1 \geq v_{m+1} = \frac{\tilde{k}_0}{m+1} \geq \frac{k_0}{m+1} = \frac{m}{m+1} v_m. \qquad (14)$$

Here there is $\mu_m \approx 0$, $v_m \approx 1$. Is not general, assume that $m$ is large enough (e.g. $m = 10^5$, At this time can use the computer software validation for all $n < 10^5$ meet the conditions of theorem 2 in $A_n$!). From formule (13), (14), we obtain

$$\mu_{m+1} = \frac{\tilde{m}_0}{m+1} \approx 0, \quad v_{m+1} = \frac{\tilde{k}_0}{m+1} \approx 1. \qquad (15)$$

Indicating that $\tilde{m}_0$ is trivial compared with $m+1$, $\tilde{k}_0$ is very close to $m+1$. So the proposition of theorem 2 is true for $n = m+1$.

Finally, according to the principle of mathematical induction, theorem 2 for any natural number $n \in N$ set up. This is the proof of theorem 2.

The conclusion of theorem 2 is very important. For any natural number $n \in N$, we can chose a large enough numbers $m \in N$ and structure matrix $A_m$ defined by(2), such that $\{f_{A_m}(j)\}_{j=2}^{m}$ is almost greater than 1, and "1, 0" are almost close to the end of $\{f_{A_m}(j)\}_{j=2}^{m}$. Thus we can infer all even numbers that less than or equal to $2n$ will be appeared in the elements of $A_m$. And so $2n$ must be the difference between the two prime numbers $p_1, p_2$, that is $2n = p_1 - p_2$.

As the application of theorem 2, we now describe and prove the following theorem.

**Theorem 4** Any even number can be represented as always in the form of at least two groups the difference between the two odd prime numbers. That is, for $n \in N$, there are prime numbers $p_1, p_2, q_1, q_2 \in P$, such that $2n = p_1 - p_2, 2n = q_1 - q_2$, and $p_1 \neq q_1$.

**Proof** For any $n \in N$, Might as well set $n$ is large enough, we structure $(4n) \times (4n)$ matrix $A_{4n}$ as fourmula (2). According to the theorem 2, the matrix master characteristic subsequence $\{f_{A_{4n}}(j)\}_{j=2}^{4n}$ of $A_{4n}$ will satisfy that it is almost greater than 1, and "1, 0" are almost close to the end of $\{f_{A_{4n}}(j)\}_{j=2}^{4n}$.

Since $n+1$ is not near the $4n$, and so $f_{A_{4n}}(n+1)$ is not near the $f_{A_{4n}}(4n)$, and then must be $f_{A_{4n}}(n+1) \geq 2$. That is to say, the number $2(n+1-1) = 2n$ must be on the $(n+1)$th line paralleled to the diagonal of the matrix $A_{4n}$ at least 2 times.

Now assume that $2n$ appear in the location of $i, j$ and $\tilde{i}, \tilde{j}$ in matrirx $A_{4n}$. Then, $i \neq \tilde{i}$, $j \neq \tilde{j}$, and

$$i < j, \quad 2i+1 \in P, \quad 2j+1 \in P, \quad 2n = 2(j-i). \tag{16}$$

$$\tilde{i} < \tilde{j}, \quad 2\tilde{i}+1 \in P, \quad 2\tilde{j}+1 \in P, \quad 2n = 2(\tilde{j}-\tilde{i}). \tag{17}$$

Let

$$p_1 = 2j+1, \quad p_2 = 2i+1, \quad q_1 = 2\tilde{j}+1, \quad q_2 = 2\tilde{i}+1, \tag{18}$$

and so $p_1, p_2, q_1, q_2 \in P$ are all odd prime such as $p_1 \neq q_1, 2n = p_1 - p_2$ and $2n = q_1 - q_2$.

This is the proof of theorem 4.

## §3 Thinking about Polignac Conjecture

What we can get from Theorem 4? Actually, theorem 4 is said all the even numbers can be said the difference between the two odd prime numbers, and the prime numbers for at least two groups. This fact makes us relate to the following two points:

First, the conclusion can be thought of as the dual form of Goldbach conjecture, the Goldbach conjecture is to prove that large even number can be expressed as the sum of two odd prime numbers. Here has proved even number can be expressed as at least two groups of the difference of the forms of two odd prime numbers (theorem 4). Is this incidental coincidence? In fact, the author in a web preprint in a similar method gives proof of the twin prime conjecture [25].

Second, this conclusion(theorem 4) answers the most weak form about the Polignac conjecture. That is for any natural number $k$, there are at least two pair prime numbers $(p, p+2k)$.

Based on these associations, we seek to a proof way to Polignac conjecture. It will engage in the work of the discussion below. Note that theorem 4, for any natural number $k$, We have at least two prime numbers $(p, p+2k)$ and $(q, q+2k)$ will be existence. Obviously, Polignac conjecture means that there are infinitely many primes $(p_i, p_i + 2k)$. Therefore, Polignac conjecture can be converted into the equivalence theorem.

**Theorem 5** Polignac conjecture is equivalent to any odd number of prime numbers $(p, q)$, there must be a natural number $k(>0)$ such as $(p+2k, q+2k)$ is a pair of primes.

**Proof** Our goal is to prove the proposition of this theorem equivalence with Polignac conjecture.

（1）Assuming Polignac conjecture be right. For nay odd prime numbers $(p, q)$, it might as well have $p \leq q$. If $p = q$, then the proposition is true obviously. In case $p < q$, let $2m = q - p$. According to Polignac conjecture, for this $m$, there are infinitely many primes

$$\{(p_i, q_i)\}_{i=1}^{\infty}, \ q_i = p_i + 2m. \ (p, q) = (p, p + 2m)$$

is also contain among them of course. Thus we have

$$2m = q - p = q_i - p_i, \ i = 1, 2, \cdots \quad (19)$$

But $q_i$ could not all smaller than $q$, so we chose $q_{i_0} > q$, then $p_{i_0} > p$ (Note: $q_{i_0} - q = p_{i_0} - p$!). Let $q_{i_0} - q = p_{i_0} - p = 2k_0$, then $k_0 > 0$ is a natural number, such as $(p + 2k_0, q + 2k_0) = (p_{i_0}, q_{i_0})$ being a pair primes. This shows that theorem 5 conclusion is true.

（2）Assuming that theorem 5 conclusion is correct. For nay natural number $m$, by Theorem 4, there are odd prime numbers $p, q$, such as

$$2m = p - q, \ p > q. \quad (20)$$

If Polignac conjecture is not true for $m$, then the pair of odd prime numbers $(p, q)$ is finite pairs in formula (20). Let there is $t_m$ pairs $\{(p_i, q_i)\}_{i=1}^{t_m}$:

$$p_i = q_i + 2m, \ p_1 < p_2 < \cdots < p_{t_m}, q_1 < q_2 < \cdots < q_{t_m}. \quad (21)$$

However, apply the theorem 5 conclusion to $(p_{t_m}, q_{t_m})$, there is a natural number $\tilde{k}(>0)$ such as $(p_{t_m} + 2\tilde{k}, q_{t_m} + 2\tilde{k})$ being a pair of primes. It is different from the outside of the primes $\{(p_i, q_i)\}_{i=1}^{t_m}$, this creates a contradiction. So Polignac conjecture should be true.

This means that Polignac conjecture is equivalence with the proposition of Theorem 5.

Theorem 5 is very important. Let $a \in N$ be the odd number, and consider the set

$$S_a = \{k \in N \mid 2k + a \in P\} \quad (1)$$

Due to the prime number set $P$ is infinite set, so $S_a$ is also infinite sets. For nay odd prime numbers $(p, q)$, Consider the corresponding sets $S_p, S_q$, it is easy to know the following conclusions.

**Theorem 6** Polignac conjecture is equivalent to any odd number of prime numbers $(p, q)$, the intersection $S_p \cap S_q$ of sets $S_p$ and $S_q$, must contain at least one natural number except 0.

**Proof** It is the direct inference of theorem 5.

It should be noted that we do not prove that Polignac conjecture, but we found the two equivalent propositions. It is easy to verify that $4 \in S_3 \cap S_5$, $5 \in S_3 \cap S_7$, $9 \in S_5 \cap S_{13}$ satisfy theorem 6 points out features. In fact, we can prove the following theorem.

**Theorem 7** Let $m \in N$, then $S_3 \cap S_{2m+3}$ must contain at least one natural number except 0.

**Proof** For $m \in N$, by the theorem 4, there are odd prime numbers $p_1, p_2, q_1, q_2 \in P$, such as $p_2 \neq q_2, 2m = p_1 - p_2$, and $2m = q_1 - q_2$. That is $p_1 = 2m + p_2, q_1 = 2m + q_2$.

Let $p_2 = 2k_0 + 3$, $q_2 = 2\tilde{k}_0 + 3$, then $k_0 \neq \tilde{k}_0$, and

$$p_1 = 2m + 2k_0 + 3 \in P, \quad p_2 = 2k_0 + 3 \in P, \quad q_1 = 2m + 2\tilde{k}_0 + 3 \in P, \quad q_2 = 2\tilde{k}_0 + 3 \in P.$$

Then according to (1), we have $k_0, \tilde{k}_0 \in S_3$, $k_0, \tilde{k}_0 \in S_{2m+3}$, and so $k_0, \tilde{k}_0 \in S_3 \cap S_{2m+3}$.

Therefore, $S_3 \cap S_{2m+3}$ must contain at least one natural number except 0. This is the proof of Theorem 7.

It must be pointed out that despite the result, so far the author can't generally prove to any odd number of prime numbers $(p, q)$, $S_p \cap S_q$ at least contains a natural number other than zero. But we can study some special cases or additional conditions as below.

**Theorem 8** For nay odd number $a \in N$, there is a number greater than 1 in $S_a \cap S_{a+4}$.

Here, the so-called "number" is greater than 1 meaning is to say, the figure relative to 1 it is very big, such as 100, 100 00, are considered to be "number greater than 1".

**Proof** Let odd number $a = 2n + 1, n \in \in N$, Using mathematical induction with $n$. First, when $n = 1$, we have

$$S_3 = \{k \in N \mid 2k + 3 \in P\} = \{1, 2, 4, 5, 7, 8, 10, 13, 14, 17, 19 \cdots\} \quad (22)$$

$$S_7 = \{k \in N \mid 2k + 7 \in P\} = \{2, 3, 5, 6, 8, 11, 12, 15, 17, 18, 20, 23, \cdots\} \quad (23)$$

And thus $2, 5, 8, 17 \in S_3 \cap S_7$, as well as $649907 \in S_3 \cap S_7$, So the proposition is true for $n = 1$.

Second, we make for inductive hypothesis: If the proposition is true for $n = k \in N$, then there is a number greater than 1 in $S_{2k+1} \cap S_{2k+5}$.

Now proves that the proposition is also founded on for $n = k + 1 \in N$, that is, whether there is a number greater than 1 in $S_{2k+3} \cap S_{2k+7}$ or not.

According to the inductive hypothesis, there is a number greater than 1 in $S_{2k+1} \cap S_{2k+5}$. So there is $k_0 \in N$ such as $k_0 \in S_{2k+1} \cap S_{2k+5}$, and $k_0$ is greater than 1. Thus

$$2k_0 + 2k + 1 \in P, \quad 2k_0 + 2k + 5 \in P \qquad (24)$$

It is noted that

$$2k_0 + 2k + 1 = 2(k_0 - 1) + 2k + 3, \quad 2k_0 + 2k + 5 = 2(k_0 - 1) + 2k + 7, \qquad (25)$$

and by (1), (24), (25), we have

$$k_0 - 1 \in S_{2k+3}, \quad k_0 - 1 \in S_{2k+7}, \quad k_0 - 1 \in S_{2k+3} \cap S_{2k+7}. \qquad (26)$$

Since $k_0$ is greater than 1, and thus $k_0 - 1$ is also greater than 1. This is a proof of that there is a number greater than 1 in $S_{2k+3} \cap S_{2k+7}$. Therefore, the proposition is also true for $n = k+1$.

According to the principle of induction, we proved that for any odd natural number $a \in N$, there is a number greater than 1 in $S_a \cap S_{a+4}$.

Theorem 8 conclusion has very important application. Because theorem 8 prove the odd $a \in N$ can freely choose, so we have

**Theorem 9** There are infinitely many groups $(p, p+4)$ such as $p$ and $p+4$ is a pair of prime numbers.

**Proof** Using reduction to absurdity. If theorem 9 is not true, then the prime formed $p$ and $p+4$ has only a limited group. Assume that $(p_0, p_0 + 4)$ be one of the largest group of prime pair $(p, p+4)$, consider the intersection of two sets $S_{p_0}, S_{p_0+4}$. By the theorem 8, there is a number greater than 1 in $S_{p_0} \cap S_{p_0+4}$. So there are $k_0 (>1) \in S_{p_0} \cap S_{p_0+4}$, such as

$$2k_0 + p_0 \in P, 2k_0 + (p_0 + 4) \in P. \qquad (27)$$

Let $q_0 = 2k_0 + p_0$, then $(q_0, q_0 + 4)$ is a new pair of prime numbers such as $(p, p+4)$. It is a pair greater than $(p_0, p_0 + 4)$ of primes formed $(p, p+4)$. This is a contradiction with $p_0, p_0 + 4$ is the biggest groups formed $(p, p+4)$. Therefore theorem is true, that is, there are infinitely many groups $(p, p+4)$ such as $p$ and $p+4$ being a pair of prime numbers.

This is the proof of theorem 9.

In this way, we prove that the primes such as $(p, p+4)$ has the infinite set. This certificate and the author is in another article to prove that the twin prime conjecture problem completely similar. In generally, we have the following theorem.

**Theorem 10** For any $m \in N$, If there be a number greater than 1 in $S_3 \cap S_{2m+3}$, then for any odd number $a \in N$, there would be a number greater than 1 in $S_a \cap S_{a+2m}$.

**Proof** Let the odd number $a = 2n+1, n \in \in N$, Using mathematical induction with $n$.

First, when $n = 1$, according to the theorem of assumptions: There is a number greater than 1 in $S_3 \cap S_{2m+3}$. So the proposition is true for $n = 1$.

Second, we make for inductive hypothesis: If the proposition is true for $n = k \in N$, then there is a number greater than 1 in $S_{2k+1} \cap S_{2k+2m+1}$.

Now proves that the proposition is also founded on for $n = k+1 \in N$, that is, whether there is a number greater than 1 in $S_{2k+3} \cap S_{2k+2m+3}$ or not.

According to the inductive hypothesis, there is a number greater than 1 in $S_{2k+1} \cap S_{2k+2m+1}$. So there is $k_0 \in N$ such as $k_0 \in S_{2k+1} \cap S_{2k+2m+1}$, and $k_0$ is greater than 1. Thus

$$2k_0 + 2k + 1 \in P, \quad 2k_0 + 2k + 2m + 1 \in P \tag{28}$$

Since

$$2k_0 + 2k + 1 = 2(k_0 - 1) + 2k + 3, \quad 2k_0 + 2k + 2m + 1 = 2(k_0 - 1) + 2k + 2m + 3, \tag{29}$$

and from the formula (1), (28), (29), we have

$$k_0 - 1 \in S_{2k+3}, \quad k_0 - 1 \in S_{2k+2m+3}, \quad k_0 - 1 \in S_{2k+3} \cap S_{2k+2m+3}. \tag{30}$$

On account of $k_0$ being greater than 1, and thus $k_0 - 1$ is also greater than 1. This is a proof of that there is a number greater than 1 in $S_{2k+3} \cap S_{2k+2m+3}$. Therefore, the proposition is also true for $n = k+1$.

According to the principle of induction, we proved that for any odd natural number $a \in N$, there is a number greater than 1 in $S_a \cap S_{a+2m}$.

Similar to theorem 9, we have

**Theorem 11** For any $m \in N$, If there be a number greater than 1 in $S_3 \cap S_{2m+3}$, then there are infinitely many groups $(p, p+2m)$ such as $p$ and $p+2m$ is a pair of prime numbers.

**Proof** Using reduction to absurdity. If theorem 11 is not true, then the prime formed $p$ and $p + 2m$ has only a limited group. Assume that $(p_0, p_0 + 2m)$ be one of the largest group of prime pair $(p, p + 2m)$, consider the intersection of two sets $S_{p_0}, S_{p_0+2m}$. By the theorem 10, there is a number greater than 1 in $S_{p_0} \cap S_{p_0+2m}$. So there exist $k_0(>1) \in S_{p_0} \cap S_{p_0+2m}$, such as

$$2k_0 + p_0 \in P, \quad 2k_0 + (p_0 + 2m) \in P, \ k_0 > 0. \tag{31}$$

Let $q_0 = 2k_0 + p_0$, then $(q_0, q_0 + 2m)$ is a new pair of prime numbers such as $(p, p+2m)$. It is a pair greater than $(p_0, p_0 + 2m)$ of primes formed $(p, p + 2m)$. This is a contradiction with $p_0, p_0 + 2m$ is the biggest groups formed $(p, p + 2m)$. Therefore theorem is true, that is, there are infinitely many groups $(p, p + 2m)$ such as $p$ and $p + 2m$ being a pair of prime numbers.

This is the proof of theorem 11. And it is different with theorem 9 in some sense, actually, this theorem is established under certain conditions.

### §4 Remarks

In summarized this paper, we mainly got the following conclusions:

（1）In the case of poor depicting two primes said even we built matrix master characteristics number and the master characteristic sequence (Definition 1), it brought convenience to prove theorem 4 in this paper.

（2）We proved that any even number can be expressed as the difference of odd prime number at least two groups (that is, Theorem 4). This is the result of comforting.

（3）We proved that there are infinitely many odd prime numbers $(p, p + 4)$ (Theorem 9), that is Polignac conjecture is true for $k = 2$.

（4）For $m \in N$, In case $S_3 \cap S_{2m+3}$ contains a number greater than 1, then there are infinitely many odd prime numbers $(p, p + 2m)$ (Theorem 11). This is Polignac conjecture is true for $k = m$ under the conditional restriction.

Here "conditional restriction" means the assumption that "for $m \in N$, if there is a number greater than 1 in $S_3 \cap S_{2m+3}$". In theorem 7, We only proves that there are at least two natural numbers in $S_3 \cap S_{2m+3}$, do not prove that it always has a number greater than 1 in it. In other words, there is a long way between our results and the Polignac conjecture. Polignac conjecture is always an unsolved problem unless someone could prove that in any case $S_3 \cap S_{2m+3}$ always contains a number greater than 1 for any $m \in N$.

But for a certain $m_0 \in N$, once the validation there is some number far greater than 1 in $S_3 \cap S_{2m_0+3}$, you can immediately use theorem 11 that Polignac conjecture in the conclusions of this $m_0$ is formed. For example, on the basis of theorem 8 and its Proof, you can take $m_0 = 2$, so theorem 9 is established.

And for $m_0 = 1$, It can verify $7743127 \in S_3 \cap S_5$ is the number of greater than 1. Therefore, in case of $m_0 = 1$, Polignac conjecture is true. This is the beginning of this article describes the famous twin prime conjecture. That is to say, this article does not ultimately proved to Polignac conjecture, but is seems to be solved the twin prime conjecture.

Anyway, if we were success in this paper, of course, it would be a very encouraging and comforting thing. The authors believe that people will eventually break Polignac conjecture.

# Chapter 3:

# Researches on the Twin Primes

**Baoshan Zhang** [①]

(School of Mathematics and Statistics, Nanjing Audit University, Nanjing, China , Zip code: 211815)

Email: bszhang@nau.edu.cn bszhang8@163.com

**Abstract**：The existence of the twin primes problem is one of the world problems in number theory. This article mainly as a result of any natural number $a \in N$, $S_a \cap S_{a+2}$ is not empty number set, and there are far more than 1 number in $S_a \cap S_{a+2}$, where

$$S_a = \{ k \in N \mid 2k + a \in P \},$$

and P is a prime number set, N is natural number set. we prove that there are an infinite number of twin prime, and then solve the problem of the twin primes in number theory.

**Keywords**　Prime numbers；Composite numbers；Twin primes；The twin prime conjecture

**Classification**：O15.156　　**MSC**：11A41

## §1 Definition and Lemmas

It is well known that natural number set $1, 2, 3, \cdots, n, \cdots$ can be classified according to different characteristics, thus get the new data set, such as the even set, collection of odd prime number set, etc. For the convenience of first presents the definition of prime number and composite number.

**Definition 1** if $n$ has no natural numbers to remove 1 and itself natural factor, $n$ is called a **prime number**; Other than 1 and a prime natural number is called a **composite number**.

According to the above definition, easy to know, $2, 3, 5, 7, 11, 13, \cdots$, are prime Numbers, and other than the two prime numbers are odd numbers, 1 is neither a prime number nor a composite. $4, 6, 8, 9, 10, 12, 14, 15, \cdots$, are composite numbers, and the even numbers except 2 are all composite numbers.

**Definition 2** if $q$, $p$ are all prime numbers, and $|p - q| \leq 2$, It says that they are **twin prime number**.

It is easy to know, $2, 3$ is the smallest prime twins, other twin primes $q$, $q+2$ are shaped like the odd prime numbers.

Usually, using $N$ express the natural number set, represented by $P$ or $\{p_k\}$ said collection of prime number, where $p_k$ is according to the natural order of the $k$-th prime Number. Here are some common results of elementary number theory in [1 ~ 2].

---

**Lemma** 1  It has an infinite number of primes, that is , Primes sequence $\{p_k\}$ is strictly monotone increasing sequence in $N$.

**Lemma** 2  Let $p$ is a prime, for any natural number $m$, or $m, p$ are relatively prime, or $m$ is divisible by $p$, namely

$$(m, p) = 1 \quad \text{or} \quad p \mid m. \tag{1}$$

**Lemma** 3（Twin Prime Conjecture [1~3]）There exist an infinite of twin primes such as $q, q+2$.

Lemma 1 ~ 2 is well known, and it has yet to be lemma 3 is proof of suspense.

On May 14, 2013, Nature magazine reported online Chinese mathematician Zhang Yitang proves that "there are an infinite number of primes less than 70 million for" poor, this research then is considered on the twin prime conjecture the ultimate number theory problem made a major breakthrough, even some people think that their impact on the community will be more than the "1 + 2" trained Chen jingrun proved [4~6]. In fact, the author has been the use spare time thinking about the twin prime conjecture in recent years, 9 ~ 11 months in 2012 found the proof of lemma 3 train of thought and simple proof is given. Based on this, this paper gives the proof of lemma 3.

## §2 The Main Results and Proofs

It $N$ express as the natural number set, $P$ as the set of prime, Let $a \in N$ be an odd natural number, Considering the number set

$$S_a = \{ k \in N \mid 2k + a \in P \} \tag{2}$$

Due to the prime number $P$ set is infinite set, so $S_a$ is also infinite sets. We can prove the following result:

**Theorem 1**  Let $a \in N$ be an odd natural number, If $k \in S_{a+2}$, then $k + 1 \in S_a$. If $2 \leq k \in S_a$, then $k - 1 \in S_{a+2}$.

**Proof**  By implication of (2), If $k \in S_{a+2}$, then $2k + (a+2) \in P$. It pay attention to that

$$2(k+1) + a = 2k + (a+2) \in P,$$

and so $k + 1 \in S_a$.

Similarly, verifiable theorem of empress half part.

The above theorem reveal the important properties of the number set defined by (2) : adjacent odd number $a, a+2$ corresponding to the number set $S_a$, $S_{a+2}$, their elements have the close relationship characterized by the theorem 1.

Thus, what is the property by the intersection $S_a \cap S_{a+2}$ of $S_a$, $S_{a+2}$? In case of $S_a \cap S_{a+2}$ is not an empty set, it exits $k \in S_a \cap S_{a+2}$ such as $k \in S_a$, $k \in S_{a+2}$, and so

$$2k + a \in P, \quad 2k + (a+2) \in P.$$

Let $q_1 = 2k + a$, $q_2 = 2k + (a+2)$, then $q_1$, $q_2$ is the twin prime number. If for any odd number $a \in N$, $S_a \cap S_{a+2}$ is not always an empty set, then the twin prime number must have infinitely many groups, which means that lemma 3 is true.

Based on these findings, we come to discuss for any odd number $a$, $b$, and the question of whether or not $S_a \cap S_b$ to the empty set. We will prove that the following theorem.

**Theorem 2** Let $a \in N$ be an any odd natural number, $S_a \cap S_{a+2}$ is not always empty, and there are far more than 1 number of $S_a \cap S_{a+2}$.

**Proof** Let $a = 2n + 1, n \in \in N$, for $n$ by using mathematical induction.

Firstly, when $n = 1$,

$$S_3 = \{ k \in N \mid 2k + 3 \in P \} = \{1, 2, 4, 5, 7, 8, 10, 13, 14, 17, 19 \cdots\} \quad (3)$$

$$S_5 = \{ k \in N \mid 2k + 5 \in P \} = \{1, 3, 4, 6, 7, 9, 12, 13, 16, 18, 19, \cdots\} \quad (4)$$

There are clearly that $1, 4, 7, 13, 19 \in S_3 \cap S_5$, So the proposition is true for $n = 1$.

Secondly, inductive hypothesis: If the depiction of the thesis is true for $n = k \in N$, that is, $S_{2k+1} \cap S_{2k+3}$ is not empty number set, and it has a number greater than 1 in $S_{2k+1} \cap S_{2k+3}$.

Now we prove that the proposition is founded on an odd number $n = k + 1 \in N$, that is, $S_{2k+3} \cap S_{2k+5}$ is not empty number set, and it has a number greater than 1 in $S_{2k+3} \cap S_{2k+5}$.

According to the induction hypothesis, the depiction of the thesis is true for $n = k \in N$, that is, $S_{2k+1} \cap S_{2k+3}$ is not empty number set, and it has a number greater than 1 in $S_{2k+1} \cap S_{2k+3}$.

So there is $k_0 \in N$ meet with $k_0 \in S_{2k+1} \cap S_{2k+3}$, where $k_0$ is a number greater than 1. And so

$$2k_0 + 2k + 1 \in P, \quad 2k_0 + 2k + 3 \in P. \quad (5)$$

Pay attention to

$$2k_0 + 2k + 1 = 2(k_0 - 1) + 2k + 3, \quad 2k_0 + 2k + 3 = 2(k_0 - 1) + 2k + 5 \quad (6)$$

Then by (2), (5), (6) and theorem 1, we have

$$k_0 - 1 \in S_{2k+3}, \quad k_0 - 1 \in S_{2k+5}, \quad k_0 - 1 \in S_{2k+3} \cap S_{2k+5}. \quad (7)$$

Because of $k_0$ being a number greater than 1, thus $k_0 - 1$ is too a number greater than 1. It's proved that $S_{2k+3} \cap S_{2k+5}$ is not empty, and $S_{2k+3} \cap S_{2k+5}$ has a number greater than 1. Therefore, the proposition to set up for $n = k+1$.

According to the principle of induction, we proved that for $a \in N$, be an any odd natural number, $S_a \cap S_{a+2}$ is not always empty, and there are far more than 1 number of $S_a \cap S_{a+2}$.

Theorem 2 the conclusion is of great importance to application. Because in theorem 2 prove that odd number $a \in N$ can be arbitrarily selected, So the lemma 3 can be thought of as direct inference of theorem 2. In fact, we can give proof of lemma 3 as below.

**Proof of lemma 3**  With the reduction to absurdity, If the conclusion in lemma 3 is not correct, that means, only a limited set of the prime twins like $p$, $p+2$. Let $p_0, p_0 + 2$ is one of the largest group of twin prime number.

Considering the collection $S_{p_0}, S_{p_0+2}$, by theorem 2, $S_{p_0} \cap S_{p_0+2}$ is not an empty number set, and it has a number greater than 1. So there are $k_0 (>1) \in S_{p_0} \cap S_{p_0+2}$, such as

$$2k_0 + p_0 \in P, \quad 2k_0 + (p_0 + 2) \in P. \tag{8}$$

Let $q_0 = 2k_0 + p_0$, then $q_0, q_0 + 2$ are twin prime numbers, They are a set of twin prime number greater than $p_0, p_0 + 2$. It is conflict with the twin prime number $p_0, p_0 + 2$ is one of the largest twin prime.

Therefore, lemma 3 is correct, that means, There must be infinitely many groups such as $p$, $p+2$, the twin prime number. This is the proof of Lemma 3.

In this way, we prove that the twin prime conjecture, solve the problem in number theory.

## §3  Remarks

This article idea for the proof of the theorem 1, 2, make out by the author among the 9 ~ 11 months in 2012. But due to busy with other interest of research, it has not been written. It is very lucky, the evening of May 17, 2013, the author in the phoenix nets see "Chinese mathematician Zhang Yitang deciphered the twin prime conjecture" [6], inspired by the news and prompted the authors to sort out their own research. Late on May 19, began to organizing and writing the first draft of this article, on May 21, the first draft is completed, on May 23, on the basis of the first draft to complete this paper.

In the first draft, the author try to prove the following proposition.

**Proposition1**  For any two odd natural numbers $a, b \in N$, $S_a \cap S_b$ is not empty number set invariably.

This proposition is promoted form of theorem 2. If the proposition 3 is proved, and it will be proved that the prime focus away from any even number of prime numbers also. Zhang Yitang given the weak form of the twin prime conjecture the amazing results that the author is convinced that the correctness of the above propositions, the author tries to give proof that seems to be not perfect, even are not able to ensure the accuracy of the certificate, so can't converse of the public to the world.

However, the theorem 2 and the proof of lemma 3, let's see the twin prime conjecture is correct, which has an infinite number of twin primes. But unfortunately, we were not constructed in the process of proof in the form of twin prime number. The author tried to construct the twin prime number column, but so far haven't see the hope of success. In this process, however, the author discovered the primes of a rule. This is the conclusion below.

**Proposition2** Each greater than 7 primes between $p$ and $2p$, there is always the twin prime number.

This conclusion is similar to Chebyshev's theorem in number theory: prime always exists between $n$ and $2n$. Due to the level, of course, cannot give evidence that it wants and cooperation to discuss the guidance of experts and scholars.